\documentclass[12pt]{article}
\usepackage[cp1251]{inputenc}
\usepackage[russian,english]{babel}
\usepackage{latexsym,amsfonts,amssymb,amsmath,amsfonts}
\usepackage{setspace}

 \newtheorem{theorem}{Theorem}[section]

 \newtheorem{remark}[theorem]{Remark}
 \newtheorem{corollary}[theorem]{Corollary}
 \newtheorem{example}[theorem]{Example}
 \numberwithin{equation}{section}

\begin{document}

%%%%%%%%%%%%%%%%%%%%%%%%%%%%%%%%%%%%%%%%%%%%%%%%%%%%%%%%%%%%%

\begin{center}

{\bf Point sources identification problems  with  pointwise overdetermination}

\medskip
%----------Author 1
S.G. Pyatkov, L.V. Neustroeva\footnote{ MSC Primary 35R30; Secondary 35R25;
35K57;

}

\bigskip

\parbox{10cm}{\small
 Abstract.  This article is  devoted to inverse  problems of
recovering point sources  in mathematical models of heat and mass transfer.
 The main attention is paid to well-posedness questions of these
 inverse problems with pointwise overdetermination conditions.  We present conditions for existence and uniqueness of
solutions to the problem, display  non-uniqueness examples, and, in model situations, we give estimates on the number of measurements that allow completely identify sources and their locations. The
results rely on asymptotic representations of Green functions of the corresponding elliptic problems with a parameter. They can be used in constructing  new numerical algorithms for determining a
solution.}

\end{center}

%%%%%%%%%%%%%%%%%%%%
%Верхние колонтитулы
%\markboth{S.G. Pyatkov}{Source identification  problems with  pointwise overdetermination}

% Заголовки разделов формируются при помощи команд  \section{}, \subsection{}, \subsubsection{}

\section*{Introduction}
\hspace{0.7 cm} Under consideration are inverse problems  with pointwise ovedetermination of
recovering  sources in mathematical models of heat and mass transfer.
We describe theoretical results  in the case of second order parabolic equation
\begin{equation}\label{1}
u_t+ A(x,D)u=f=\sum_{i=1}^r f_{i}(t,x)q_i(t)+f_0,\ (t,x)\in Q=(0,T)\times G\  (G\subset {\mathbb R}^n).
 \end{equation}
 The operator $A$ is an elliptic operator  representable as
  \begin{equation*}
A(x,D)u=-\sum_{i,j=1}^{n}a_{ij}(x)u_{x_ix_j} +\sum_{i=1}^n a_i(x)u_{x_i}+a_0(x)u.
 \end{equation*}
 The equation \eqref{1} is furnished with the initial and boundary conditions
\begin{equation}\label{2}
u|_{t=0}=u_0,\ \ Bu|_{S}= g(t,x),\  S=(0,T)\times \Gamma, \ \Gamma=\partial G,
 \end{equation}
 where $Bu=u$ or   $Bu=\frac{\partial u}{\partial N}+\sigma u=
 \sum_{i,j=1}^{n} a_{ij}(t,x)u_{x_j}(t,x)\nu_i+\sigma(t,x)u(t,x)$,
where $\vec{\nu}=(\nu_1,\ldots,\nu_n)$ is the outward unit normal to $\Gamma$.
The unknowns in \eqref{1}, \eqref{2} are a solution  $u$ and the function $q_i(t)$ ($i=1,2,\ldots, r$) occurring into
the right-hand side of \eqref{1}. The overdetermination conditions for recovering the functions $\{q_i\}_{i=1}^r$  are as follows:
\begin{equation}\label{3}
u|_{x=b_i}=\psi_i(t),\ \ i=1,2,\ldots,s,
 \end{equation}
where  $\{b_i\}$ is a collection of points in $G$ or on $\Gamma$.

 These problems   arise in  mathematical modelling of
heat and mass transfer processes, diffusion, filtration, and in many other fields (see \cite{mar,ozi,ali}).  A function $f$ on the right-hand side of \eqref{1}
  is referred to as a source function.
Generally,  there are two   types of inverse source problems, one of them is the problem   of recovering
point sources in which  $f_{i}(t,x)=\delta(x-x_i)$, with $\delta$ the Dirac delta function. In the latter case
the functions  $f_{i}(t,x)$ are some smooth or at least summable functions. In this case the function $f$ can be called the distributed source function.
First of all, we should refer to the fundamental articles by
Prilepko A.I. and his followers. In particular, an existence
and uniqueness theorem for solutions to the problem of
recovering the source  $f(t,x)q(t)$ with the overdetermination
condition $u(x_{0},t)=\psi(t)$ $(x_{0}$ is a point in  $G$) is
established in \cite{pri,sol}. The H\"{o}lder spaces serve as the basic spaces  in these
articles. The results were generalized in the book
\cite[\S\,6.6,\S\,9.4]{pri2}, where the existence theory for
the problems \eqref{1}-\eqref{3} was developed in an
abstract form with the operator  $A$ replaced with $-L$, $L$ is
generator of an analytic semigroup. The main results employ the
assumptions that the domain of $L$ is independent of time and
the unknown coefficients occur into the lower  part of the
equation nonlinearly. Under certain conditions,  existence and
uniqueness theorems were proven locally in time in the spaces
of  functions continuously differentiable with respect to time.
    Many results in
the case of  $n=1$ are exhibited in  \cite{bib3}.  The problems \eqref{1}-\eqref{3} were considered in author's articles in \cite{pya3}-\cite{pya7},  where conditions on the data were weakened in
contrast to those in \cite[\S\,9.4]{pri2} and the solvability questions were treated in the Sobolev spaces. Note that in the article \cite{pya7} the unknowns occur in the right-hand side on \eqref{1}
nonlinearly. In this case the right-hand side is of the form $f(t,x,u,\nabla u, \vec{q}(t))$ (see also  \cite{pri}).
 Proceed with the point sources problems, which, in contrast to the distributed sources case,
 is always ill-posed. In this case in the problem \eqref{1}-\eqref{3} with point sources, i.~e., $f=\sum_{i=1}^r \delta(x-x_i)q_i(t)+f_0$,  the intensities $q_{i}(t)$
of point sources, their  locations $x_{i}$ and the number $m$
are quantities to be determined. There is a small number of theoretical results devoted to the
 solving these inverse problems.  The
main results are connected with numerical methods of solving the problem and many of them are far from justified. The problem is ill-posed and examples when the problem is not solvable or has many
solutions are easily constructed. Very often the methods rely on reducing the problem to an optimal control problem and minimization of the corresponding objective functional \cite{16,18,pen,ozi}.
However, it is possible that the corresponding functionals can have many local minima. Some theoretical results devoted to the problem \eqref{1}-\eqref{3} are available in  \cite{4}-\cite{11}. The
stationary case is treated in \cite{bad}, where the Dirichlet data  are complemented with the Neumann data  and these data allow to solve the problem on recovering the number of sources, their
locations, and intensities  using test functions and a Prony-type algorithm. Similar results are obtained in the multidimansional case for parabolic source problem \cite{17} and thereby the
identifiability of point sources is proven in the case of Cauchy data on the boundary of a spatial domain (i.~e., the Dirichlet data on the boundary in addition to the Neumann data  are given).
 The model problem
\eqref{1}-\eqref{3} ($G={ \mathbb R}^{n}$) is considered in
\cite{11}, where the explicit representation of solutions to
the direct problem (the Poisson formula) and an auxiliary
variational problem are employed to determine numerically  the
quantities $\sum_{i}q_{i}r_{ij}^{l}$ (here
 $q_{i}(t)=const$ for all  $i$ and
$r_{ij}=|x_{i}-b_{j}|$, $l=1,2,\ldots$). The quantities found  allow to determine the points $\{x_{i}\}$ and intensities $q_{i}$ (see Theorem 2 and the corresponding algorithm in  \cite{11}).  In the
one-dimensional case, uniqueness theorem for solutions to the problem \eqref{1}-\eqref{3} with $m=1$ is stated in \cite{4}. Similar results are presented also in \cite{pya1}. To define  a solution
uniquely (intensities and source locations) in the one-dimensional case, we need the condition that the sources and measurement points $\{b_i\}$ alternate and this requires some unavailable
information in a practical situation. Non-uniqueness examples  in the problems of recovering of point sources are presented in \cite{neus}. Some numerical methods   used in solving the problems of
recovering point sources  are described in \cite{zhou}-\cite{pek}. The bibliography in the stationary case can be found in \cite{hon}.

This article is actually a survey of the results obtained in the articles \cite{pya1}, \cite{pya111}-\cite{pya112}. In some model cases, we discuss
 the uniqueness and existence results for the problem \eqref{1}-\eqref{3}. The results uses a new approach based on asymptotic representations of
 the Green function of the corresponding elliptic problem.

   \section{Preliminaries}
\hspace{0.7 cm}

Let  $E$ be a~Banach space. The~symbol  $L_p(G;E)$ ($G$
   is a~domain in ${\mathbb  R}^n$) stands for the~space of  measurable functions defined on  $G$ with values in
    $E$ and a~finite norm  $\|\|u(x)\|_{E}\|_{L_p(G)}$.
The notations of the~Sobolev spaces  $W_p^s(G;E)$ and $W_p^s(Q;E)$ are conventional (see \cite{ama}). If $E={\mathbb  R}$ or
    $E={\Bbb  R}^n$ then the~latter space is denoted by  $W_{p}^{s}(Q)$. The~definitions of the~H\"{o}lder spaces
    $C^{\alpha,\beta}(\overline{Q}),
    C^{\alpha,\beta}(\overline{S})$ can be found, for example, in
   \cite{lad}.
By the~norm of a~vector, we mean the~sum of the~norms of its coordinates.
Given an interval
$J=(0,T)$, put
    $W_p^{s,r}(Q)=W_p^{s}(J;L_p(G))\cap L_p\bigl(J;W_p^r(G)\bigr) $ and
  $W_p^{s,r}(S)=W_p^{s}(J;L_p(\Gamma))\cap
    L_p\bigl(J;W_p^r(\Gamma)\bigr)$.
Let     $(u,v)=\int\limits_{G}u(x)v(x)\,dx$ and denote by
  $B_{\delta}(b)$ the~ball of radius  $\delta$ centered at   $b$.
 The symbol $\rho(X,Y)$ stands for the distance between the sets $X,Y$.

   \section{The one-dimensional case}

We consider the question of recovering the source function of a
special form in the one-dimensional
advection–dispersion–reaction  equation
\begin{equation} \label{eq1}
\begin{array}{c}
Lu=u_t-L_{0}u=\sum_{i=1}^{r}q_{i}(t)\delta(x- x_{i})+f_0(x,t)=f,\ \vspace{4pt} \\
(x,t)\in G\times (0,T)=Q,\  G=(a,b),
\end{array}
\end{equation}
 where  $ L_{0}u=a_2(x)u_{xx}-a_1(x)u_{x}-a_0(x)u$,   $\delta$ is the Dirac delta-function and $-\infty\leq a<b\leq \infty$, $T\leq \infty$.
Here   $u(x,t)$ is the pollutant  concentration and $q_{i}(t)$
are the intensities of point sources with coordinates $x_{i}$
($i=1,2,\ldots,m$)  (see the description of models, for
instance, in \cite{mar}). In the general case the functions $q_{i}(t)$, the points $x_i$, and the number $r$ as well as a
solution $u$ are assumed to be unknown. The equation
(\ref{eq1}) is furnished with the initial and boundary
conditions
\begin{equation} \label{eq2}
B_{j}u=\varphi_{j}(t),\ \  j=1,2,\  u(x,0)=u_{0}(x),
\end{equation}
where, in the case of   $a\neq -\infty$,  $B_{1}u=u(t,a)$ or
$B_{1}u=u_{x}+\sigma_1 u|_{x=a}$, respectively, in the case of  $b\neq +\infty$,
$B_{2}u=u(t,b)$ or $B_{2}u=u_{x}+\sigma_2 u|_{x=b}$ ($\sigma_i=const$). In the case of
 $a=-\infty$ or $b=+\infty$ the corresponding boundary condition is replaced with the condition  $\lim_{x\to -\infty}u(x,t)=0$ or $\lim_{x\to
+\infty}u(x,t)=0$, respectively, which are understood in an integral sense.
We use  the  pointwise overdetermination
conditions
\begin{equation} \label{eq3}
u(b_{i},t)=\psi_{i}(t),\ \  i=1,2,\ldots, s.
\end{equation}
Thus, the problem can be stated as follows: given functions $\psi_{i}$ $(i=1,\ldots, s)$, $u_{0}$, $\varphi_{i}$ ($i=1,2$), find a solution
 $u$ to the equation (\ref{eq1}), the functions  $q_{i}(t)$ $(i=1,2,\ldots,r)$, and  points $x_{i}$ $(i=1,2,\ldots,r)$ such that the
 equalities  \eqref{eq1}, \eqref{eq2}, and \eqref{eq3} hold.
Describe the class in which we look for a solution to the problem.
     By a solution to the problem \eqref{eq1}-\eqref{eq2}, we mean a generalized solution. For instance, if $T<\infty$ and the conditions \eqref{eq2} are the Robin boundary conditions then  a function  $u\in  L_2\big(0,T;W_2^1(a,b)\big)$ is a generalized solution to the problem \eqref{eq1}, \eqref{eq2} whenever
\begin {multline*}
\int_Q -u \psi_t +a_2(x)u_x \psi_x +(a_1+a_{2x})u_x \psi +a_0u\psi\,dxdt -\int_G u_0 \psi(0,x)\,dx+ \\
\int_0^T a_2(b)\sigma_1u(t,b)\psi)(t,b)-a_2(a)\sigma_2u(t,b)\psi)(t,b)\,dt + \\ \int_0^T \varphi_1(t)\psi(t,a)-\varphi_2(t)\psi(t,b)\,dt  =\int_Q f(t,x)\psi(t,x)\,dxdt\
\end {multline*}
for all $\psi(t,x)\in W_2^1(Q)$ such that  $supp\,\psi\in [0,T)\times G$. In other cases the definition is quite similar (see the definitions in \cite{lad}).

The problem in this statement is not correct in the Hadamard sense and we have no solvability  in classes of finite smoothness
as well as uniqueness of solutions. The examples are easy to provide.

To obtain uniqueness, we need to impose the condition that the points $\{x_i\}$ and $\{b_i\}$ alternate in a certain sense even in the case we consider a simpler problem in which the location of points $\{x_i\}$ is assumed to be known and we need to determine the intensities $\{q_i(t)\}$ and a solution $u$.

We enumerate the points in~$S_x=\{x_i\}_{i=1}^r$ and $S_y=\{b_j\}_{j=1}^s$ in increasing order and  let $a_2(x)\in W_\infty^1(a,b)$, $a_i\in L_\infty(a,b)$ ($i=0,1$),
$M_{1}\leq a_2(x)\leq M_{2}$ $\forall x\in (a,b)$ and some positive constants $M_1,M_2$.
In  this case the following theorem is valid (see \cite[Theorem 4]{pya111}).

\begin{theorem}\label{th1} Assume that one of the following condition holds;
    \begin{enumerate}%[\quad\upshape (a)]
    \item $a<x_1<x_2$, $b_i>x_2$, $i=1,2,\dots,s$,
    \item $x_{r-1}<x_r<b$, $b_i<x_{r-1}$, $i=1,2,\dots,s$,
    \item $[x_{r_1},x_{r_1+2}]\cap S_y=\emptyset$ for some  $r_1$.
    \end{enumerate}
    Then a solution to the problem  {\rm \eqref{eq1}-\eqref{eq3}} on recovering intensities
    $q_i(t)$ $(i=1,2,\dots,r)$ and a function ~$u$ is not unique.
 \end{theorem}

\begin{remark} If the points $\{x_i\}$ are known then we can arrange the points $\{b_i\}$  to ensure uniqueness of solutions to our problem.
Otherwise,  uniqueness of solution to the problem cannot be guaranteed.
Examples show that the problems of recovering the intensities $\{q_i\}$ or the problem of recovering the points $\{x_i\}$ and the intensities $\{q_i\}$
can be reduced to a collection of problems \eqref{eq1}-\eqref{eq3} with $r=1$ if the points $\{b_i\}$ are arranged properly.
However, even in this case the problem is not simple.
\end{remark}

Describe some simplest results in the case of $r=1$.
 We assume that
\begin{equation}\label{eq4}
    a_2\in C^1\big([a,b]
            \big)\cap W_1^2(a,b),\ \ a_1\in C
            \big([a,b]
            \big)\cap W_1^1(a,b),\ \ a_0\in L_{\infty}(a,b)
  \end{equation}
   in the case of the finite interval  $(a,b)$ and
          \begin{equation}\label{eq5}
    a_2\in C^1\big([a,b]
            \big)\cap W_\infty^2(a,b),\ \ a_1\in C
            \big([a,b]
            \big)\cap W_\infty^1(a,b),\ \ a_0\in L_{\infty}(a,b)
 \end{equation}
if the interval $(a,b)$ is infinite.

Write out our conditions for the data. Fix $\lambda\geq 0$.

(A) if the $i$-th condition in (\ref{eq2}) ($i=1,2$) is the
Dirichlet  condition (respectively, the Neumann condition)  then $e^{-\lambda t} \varphi_{i}\in
W_{2}^{3/4}(0,T)$ (respectively, $e^{-\lambda t}\varphi_{i}\in
W_{2}^{1/4}(0,T)$); $u_{0}\in W_{2}^{1}(a,b)$, $e^{-\lambda t} f_0\in L_{2}(Q)$.

The consistency conditions are as follows:

(B) if the boundary condition at $x=a$ ($x=b$) is the Dirichlet
boundary condition then $\varphi_{1}(0)=u_{0}(a,t)$
(respectively, $\varphi_{2}(0)=u_{0}(b,t)$).

(C)  $\psi_{i}(0)=u_{0}(b_{i})$, $i=1,2,\ldots,s$.

Conventional results imply that

\begin{theorem}\label{th2} Let  the conditions {\rm \eqref{eq4}, \eqref{eq5}} hold. If  $T=\infty$
then  there exists a parameter $\lambda _{0}\geq 0$ such that if the conditions  {\rm (A), (B)} hold
for some $\lambda \geq \lambda_{0}$ then there exists a unique solution $\Phi$
to the problem {\rm \eqref{eq1}, \eqref{eq2}}, where $q_i(t)\equiv 0$ for all $i$ such that  $e^{-\lambda t}\Phi\in
L_{2}(0,T;W_{2}^{2}(G))$, $e^{-\lambda t}\Phi_{t}\in L_{2}(Q)$.
If   $T<\infty$ and the conditions   {\rm (A), (B)} with $\lambda =0$ hold then
 there exists a unique solution $\Phi$
to the problem {\rm \eqref{eq1}, \eqref{eq2}}, where $q_i(t)\equiv 0$ for all $i$ such that
 $\Phi\in L_{2}(0,T;W_{2}^{2}(G))$, $\Phi_{t}\in L_{2}(Q)$.
\end{theorem}

This  result follows from Theorem 10.1 in \cite{agr} (see also \cite[Theorems 5.7,7.11,8.2]{denk}).

 Let the
space $\tilde{W}_{2}^{1}(a,b)$ comprise the functions in
$W_{2}^{1}(G)$ satisfying those homogeneous boundary conditions
in (\ref{eq2}) that have a sense. Thus, if both conditions in
(\ref{eq2}) are the Neumann conditions then
$\tilde{W}_{2}^{1}(a,b)=W_{2}^{1}(a,b)$. Assume that
$\tilde{W}_{2}^{-1}(a,b)$ is dual space to
$\tilde{W}_{2}^{1}(a,b)$ relatively to the inner product in
$L_{2}(a,b)$. It can be defined as the completion of
$L_{2}(a,b)$ with respect to the norm
$\|u\|_{\tilde{W}_{2}^{-1}(a,b)}=\sup_{v\in
\tilde{W}_{2}^{1}(a,b)} |(u,v)|/\|v\|_{W_{2}^{1}(a,b)}$.

 Using Theorems \ref{th2}, we can construct a function  $\Phi$ and  make the change of variables  $u=\omega+\Phi$
in (\ref{eq1})-(\ref{eq3}).
 Our inverse  problem for $r=1$
is reduced to the problem
\begin{equation} \label{eq6}
\begin{array}{c}
L\omega=\omega_t-L_{0}\omega=q_1(t)\delta(x- x_{1}),
\end{array}
\end{equation}
\begin{equation} \label{eq7}
B_{j}\omega=0,\ \  j=1,2,\  \omega(x,0)=0,
\end{equation}
\begin{equation} \label{eq8}
\omega(b_{i},t)=\tilde{\psi}_{i}(t)=\psi_{j}(t)-\Phi(b_{j},t) ,\ \  i=1,2,\ldots, s.
\end{equation}

We can state that the direct problem (\ref{eq6}), (\ref{eq7})
is solvable in the sense that (we take $T\leq
\infty$)

\begin{theorem}\label{th3} There exists a parameter $\lambda_0\geq 0$ such that if  $\lambda\geq \lambda_0 $,  $e^{-\lambda t} q_1(t)\in L_{2}(0,T)$, and the conditions {\rm (\ref{eq4}),  (\ref{eq5})} hold, then there
exists a unique solution  $\omega$ to the problem  {\rm
(\ref{eq6}), (\ref{eq7})} from the class  $e^{-\lambda
t}\omega\in L_{2}(0,T; W_{2}^{1}(a,b))$, $e^{-\lambda t}
\omega_{t}\in L_{2}(0,T; \tilde{W}_{2}^{-1}(a,b))$. If
$q_1(t)\geq 0$ on $(0,T)$ then a solution $\omega$ is nonnegative
in $Q$. Moreover, $e^{-\lambda t}\omega\in
L_{2}(0,T;W_{2}^{2}(G_{\varepsilon}))$, $e^{-\lambda
t}\omega_{t}\in L_{2}(0,T; L_{2}(G_{\varepsilon}))$ for all
$\varepsilon>0$ $(G_{\varepsilon}=\{x\in G:\ |x-x_i|\geq \varepsilon \ \forall i\})$.
\end{theorem}

This theorem  can be proven with the use of \cite[Theorem 14.2,
Theorem 17.1]{ama} (see also \cite{lad}).

Denote $Q_{\varepsilon}=(0,T)\times G_{\varepsilon}$ and   $r(\xi)=1\big/\sqrt{a_2(\xi)}$.
The following theorem states that the problem  \eqref{eq1}-\eqref{eq3}, where
$r=1$ and $s=1$ and we determine a solution $u$ and the intensity $q_{1}$ is uniquely solvable.
As before, we construct an auxiliary function $\Phi$ and make the change of variables
  $\omega=u-\Phi$.

 Consider the function $V_{\delta}(t)$ ($\delta>0$) such that
    $$
    {\cal L}(V_{\delta})(\lambda)=e^{-\sqrt{\lambda}\delta}/\sqrt{\lambda},
    \ \  \lambda^{\alpha}=|\lambda|^{\alpha}e^{i\alpha\arg \lambda}.
    \ -\pi<\arg \lambda<\pi,
    \ \lambda=\sigma+i\gamma,
    $$
    where ${\cal L}$ is the Laplace transform.
         Define the class of functions
    $$
    H_{\delta}=\big\{\psi(t)=\int_0^t\psi_0(\tau)V_{\delta}(t-\tau)\,d\tau:
                      \ \psi_0\in L_2(0,T)
               \big\}.
    $$

\begin{theorem}\label{th4} Let   $r=1$,  $s=1$, $\tilde{\psi}_1\in H_{\delta_0}$ with $\delta_0=\big|\int_{x_1}^{b_1}r(\xi)\,d\xi\big|$,
the conditions  {\rm \eqref{eq4},  \eqref{eq5}}, {\rm (A)-(C)}, where $\lambda=0$ and    $T<\infty$.
 Then there exists a unique solution ~$(u,q_1)$ to the problem  \eqref{eq1}-\eqref{eq3} such that
    $u\in L_2\big(0,T;W_2^1(a,b)\big)$,
    $(u-\Phi)_{t}\in L_2\big(0,T;\tilde{W}_2^{-1}(a,b)\big)$,
        $u\in W_2^{1,2}(Q_{\varepsilon})$ for all     $\varepsilon>0$.
Assume that   $T=\infty$,  $r=1$,  $s=1$ and the conditions  {\rm \eqref{eq4},  \eqref{eq5}} hold. Then there exists a parameter
$\lambda_{0}>0$ such that if $\lambda\geq \lambda_0$, the conditions  {\rm (A)-(C)} hold, and the function  $\tilde{\psi}_1$ admits the representations
    \begin{equation}\label{eq9}
    \tilde{\psi}_1=\int\limits_0^t\psi_0(\tau)V_{\delta_0}(t-\tau)\,d\tau,
    \  \delta_0=\big|\int\limits_{x_1}^{b_1}r(\xi)\,d\xi\big|,\  \psi_0(\tau)e^{-\lambda \tau}\in L_2(0,\infty),
    \end{equation}
    then there exists a unique solution ~$(u,q_1)$ to the problem  \eqref{eq1}-\eqref{eq3}, where  $r=1$ and $s=1$, such that
    $ue^{-\lambda t}\in L_2\big(0,T;W_2^1(a,b)\big)$,
    $(u-\Phi)_{t}e^{-\lambda t}\in L_2\big(0,T;\tilde{W}_2^{-1}(a,b)\big)$,
        $ue^{-\lambda t}\in W_2^{1,2}(Q_{\varepsilon})$ for all
        $\varepsilon>0$.
\end{theorem}

The proof can be found in \cite{neus3}.
Outline the proof. Applying the Laplace transform to the equation \eqref{eq6}, we arrive at the problem
\begin{equation} \label{eq10}
\lambda \hat{\omega}-L_{0}\hat{\omega}=\hat{q}_1(\lambda)\delta(x- x_{1}),
\end{equation}
\begin{equation} \label{eq11}
B_{j}\hat{\omega}=0,\ \  j=1,2,
\end{equation}
\begin{equation} \label{eq12}
\hat{\omega}(b_{1},t)={\cal L}(\tilde{\psi}_{1})(\lambda).
\end{equation}
From the equalities  \eqref{eq10}- \eqref{eq12}, it follows that
\begin{equation} \label{eq13}
{\cal L}(\tilde{\psi}_{1})(\lambda) = \hat{q}_1(\lambda) (\lambda -L_{0})^{-1}\delta(x- x_{1}).
\end{equation}
Next, we employ the asymptotic representation of  the function $v_0=(\lambda -L_{0})^{-1}\delta(x- x_{1})$ (see \cite{neus3}) and the properties of the Laplace transform. The asymptotic representation is of the form.
 Fix  $\delta_0\in(0,\pi)$. There exists  $\lambda_0\geq 0$ such that for all  $\lambda\in {\mathbb C}$ such that
   $|\arg(\lambda-\lambda_0)|\leq\pi-\delta_0$ the function $v_0=(\lambda -L_{0})^{-1}\delta(x- x_{1})$ is defined,
    $v\in W^1_2(a,b)$,  and on any compact set  $[c,d]\subset (a,b)$ admits the representation
     \begin{equation}\label{eq14}
        v_0(b_1)=\frac{-1}{2\sqrt{\lambda a(x_{1})}}
    \exp\Big(\!{-}\sqrt{\lambda}\,\Big|\int\limits_{x_1}^{b_1}\!\!r(\xi)\,d\xi
                                 \Big|{+}\!\!
              \int\limits_{x_1}^{b_1}\!\!r_1(\xi)\,d\xi
        \Big)\Big(1{+}O\Big(\frac1{\sqrt{\lambda}}
                         \Big)               \Big),
\end{equation}
where $r_1(\xi)= \frac{-1}2(a_2r'r-a_1r^2)(\xi)$.

 Next, we consider the problem  \eqref{eq1}-\eqref{eq3} with $r=1$ of simultaneous determination of the intensity $q_1(t)$ and
 the point $x_1$.
 We consider the case of  $r=1, s=2$ and the problem of recovering the quantities  $(u,q_1,x_{1})$.
Constructing the function  $\Phi$, we reduce the problem to the problem \eqref{eq6}-\eqref{eq8}. It is possible that
 $b_{1}=a$ (if  $a\neq -\infty)$) and, respectively, that $b_{2}=b$ (if  $b\neq +\infty)$. In the former case we take
 $B_{1}u=u_{x}(a,t)$ and in the latter
$B_{2}u=u_{x}(b,t)$.
First of all, we can say that
\begin{theorem}\label{th4}  Under the conditions  {\rm \eqref{eq4},  \eqref{eq5}}, a solution to the problem   \eqref{eq1}-\eqref{eq3}
where $r=1$ and $s=2$, of defining  $x_1$,  $u$, and   $q_1$~is unique if it is known that $x_1\in (b_1,b_2).$
  \end{theorem}

Let
$\Phi_i(\lambda)={\cal L}(\tilde{\psi}_i)$, $i=1,2$.

\begin{theorem}\label{th5} Assume that the conditions of Theorem  {\rm \ref{th3}} hold  for some $\lambda \geq \lambda_0$ {\rm (}the parameter $\lambda_0$ is that of Theorem {\rm \ref{th3}}{\rm )} and     $(u,q_1,x_{1})$~is a solution to the problem  {\rm \eqref{eq1}-\eqref{eq3}}, where $T=\infty$, from the class indicated in
Theorem {\rm \ref{th3}},    $q_1\not\equiv 0$, and
    $b_{1}<x_{1}<b_{2}$.
    Then there exists ~$\lambda_1>0$ such that on the segment
    ~$[\lambda_2,\infty)$ the sets of zeros of the functions ~$\Phi_1(\lambda)$ and
    $\Phi_2(\lambda)$ coincide and have no finite limit points, there exists the finite limit
        $$
    \lim_{\lambda\to\infty}
    \frac1{2\sqrt{\lambda}}\ln
    \frac{\Phi_1(\lambda)}
         {\Phi_2(\lambda)}=A,
        $$
    and we have the following asymptotic equalities:
     if $b_{1}= a$ and $b_{2}< b$ then
\begin{equation}\label{eq15}
    \int\limits_a^{x_1}r(\xi)\,d\xi =\frac12
                   \int\limits_a^{b_2}r(\xi)\,d\xi
 -    \frac1{2\sqrt{\lambda}}\ln
    \frac{\Phi_1(\lambda)}
         {2\Phi_2(\lambda)}-
    \frac1{2\sqrt{\lambda}}
    \int\limits_{a}^{b_2}r_1(\xi)\,d\xi+O
    \bigg(\frac1{\lambda}
    \bigg);
\end{equation}
        if  $b_{1}> a$ and $b_{2}= b$ then
\begin{equation}\label{eq16}
    \int\limits_{x_1}^{b}r(\xi)\,d\xi
 =\frac12\int\limits_{b_1}^b r(\xi)\,d\xi+\frac{1}{2\sqrt{\lambda}}
                   \int\limits_{b_1}^{b}r_1(\xi)\,d\xi
                         +    \frac1{2\sqrt{\lambda}}\ln
    \frac{2\Phi_1(\lambda)}
         {\Phi_2(\lambda)}+O\bigg(\frac1{\lambda}
    \bigg);
\end{equation}
  if $b_i\in (a,b)$ then
\begin{equation}\label{eq17}
    \int\limits_{b_1}^{x_1}r(\xi)\,d\xi
 =\frac12\int\limits_{b_1}^{b_2}r(\xi)\,d\xi-
                  \frac{1}{2\sqrt{\lambda}} \int\limits_{b_1}^{b_2}r_1(\xi)\,d\xi
                 -\frac{1}{2\sqrt{\lambda}}\ln
    \frac{\Phi_1(\lambda)}{\Phi_2(\lambda)}+O\bigg(\frac1{\lambda}\bigg).
\end{equation}
\end{theorem}

The proof can be found in \cite{neus3}.
An analog of Theorem \ref{th5} also holds in the case of a finite segment  $[0,T]$ but the quantity
 $O \bigg(\frac1{\lambda}    \bigg)$ is replaced with
 $o\bigg(\frac1{\sqrt{\lambda}}\bigg)$ (see Theorem 3 in \cite{pya1}).
The proofs rely on properties  of the Laplace transform and an asymptotic representation of the form \eqref{eq14}.

 \begin{corollary}\label{col1}
Passing to the limit as  $\lambda\to \infty$, we obtain the equality
 \begin{equation}\label{eq18}
    A= \frac{1}{2}\int_{b_1}^{b_2}r(\xi)\,d\xi - \int_{b_1}^{x_1}r(\xi)\,d\xi,
    \end{equation}
which implies that
 \begin{equation}\label{eq19}
    |A|< \frac{1}{2}\int_{b_1}^{b_2}r(\xi)\,d\xi.
    \end{equation}
We can see from \eqref{eq18} that the quantity  $x_1\in (b_1,b_2)$ is defined uniquely.
    \end{corollary}

\begin{remark} \label{rem000} The formulas  \eqref{eq15}-\eqref{eq17} can be considered as equations for
determining the point $x_1$ and can be used in numerical calculations. If $a_2(x)=const$ and $b_i\in (a,b)$ then the equation  \eqref{eq17} can be rewritten as follows:
\begin{equation}\label{eq20}
    x_1
 =\frac{b_1+b_2}{2} -
                  \frac{\sqrt{a_2}}{2\sqrt{\lambda}} \int_{b_1}^{b_2}r_1(\xi)\,d\xi
                 -\frac{\sqrt{a_2}}{2\sqrt{\lambda}}\ln
    \frac{\Phi_1(\lambda)}{\Phi_2(\lambda)}+O\bigg(\frac1{\lambda}\bigg).
\end{equation}
Numerical experiments show that the determination of $x_1$ from this equality it is quite possible  \cite{pya1} if we choose an appropriate parameter
$\lambda$.  The problem is that we need to calculate the integrals  $\Phi_i(\lambda)$. To calculate them with a good accuracy,
we should require  the inequality $\lambda\leq c/\tau$, where $\tau$ is the time step and $c$ is some small constant.
So  the increase of the parameter  $\lambda$ cannot  give more exact determination of
 $x_1$ due to errors of calculations.
\end{remark}

 \section{The multi-dimensional case}
  In this section we expose some results concerning with recovering of point sources in the multidimensional case.
  We consider a simple parabolic  equation
  \begin{equation}\label{eq21}
u_{t}+Lu = \sum_{i=1}^r q_i(t)\delta(x-x_i)+f_{0}(t,x),\ \ Lu=-\Delta u + \sum_{i=1}^{n}a_{i}(x)u_{x_{i}}+a_{0}(x)u,
 \end{equation}
 where $(x,t)\in Q=(0,T)\times G$, $G$ is a domain in ${\mathbb R}^{n}$ ($n=2,3$) with boundary $\Gamma\in C^{2}$.
 We consider three cases $G={\mathbb R}^n$, or $G={\mathbb R}_+^n=\{x:x_n>0\}$, or $G$
 is a domain with compact boundary.
The unknowns are the functions  $q_{i}(t)$. The equation
\eqref{eq21} is furnished with the initial and boundary
conditions
\begin{equation}\label{eq22}
    Bu|_{S}=g,\  \ u|_{t=0}=u_{0}(x),\ \ S=(0,T)\times \Gamma,
    \end{equation}
where either  $Bu=\frac{\partial u}{\partial \nu}+\sigma u$, or
$Bu=u$ ($\nu$ is the outward unit normal to $\Gamma$), and the
overdetermination conditions are as follows:
\begin{equation}\label{eq23}
    \ u(b_{j},t)=\psi_{j}(t),\  j=1,2,\ldots,s.
    \end{equation}
The coefficients in \eqref{eq22} are real-valued.

First, we describe our conditions on the data and present the simplest existence theorem.
 Let $\vec{a}=(a_{1},a_{2})$ for $n=2$
and $\vec{a}=(a_{1},a_{2},a_{3})$ for $n=3$. The brackets
$(\cdot,\cdot)$ denote the inner product in ${\mathbb R}^{n}$.
The coefficients in  \eqref{eq1}  are assumed to be real-valued
and
    \begin{equation}\label{q13}
a_{i}\in W_{\infty}^{2}(G)\ (i=1,\ldots,n),\  a_{0}\in L_{\infty}(G), \ \sigma\in C^{1}(\Gamma),
    \end{equation}
Fix a parameter $\lambda\in {\mathbb R}$ and assume that
\begin{equation}\label{eq603}
e^{-\lambda t}g\in W_{2}^{1/4,1/2}(S),\  \textrm{if } Bu=\frac{\partial u}{\partial \nu}+\sigma u\ (\sigma\in C^{1}(\Gamma)),
 \end{equation}
 \begin{equation}\label{eq604}
e^{-\lambda t} g\in W_{2}^{3/4,3/2}(S),\  \textrm{if } Bu=u,\ \ f_{0}e^{-\lambda t}\in L_{2}(G),\ \
 \end{equation}
 \begin{equation}\label{eq601}
 u_{0}(x)\in W_{2}^{1}(G),\ \
u_{0}(x)|_{\Gamma}=g(x,0)\  \textrm{if } Bu = u.
 \end{equation}

\begin{theorem} \hspace{-5pt}.
\label{th6} Assume that  $T=\infty$, $a_i\in L_\infty(G)$ $(i=0,1,\ldots,n)$, and
the conditions  {\rm \eqref{eq601}} hold.
 Then there exists  $ \lambda_{0}\geq 0$ such that if
 $\lambda\geq \lambda_{0}$ and the conditions
{\eqref{eq603}, \eqref{eq604}} are fulfilled then there exists a unique solution
to the problem
 \begin{equation}\label{eq602}
w_{0t}+Lw_0 =f_{0}(t,x), \ \   Bw_0|_{S}=g,\  \ w_0|_{t=0}=u_{0}(x),\
 \end{equation}
such that  $e^{-\lambda t} w_{0}\in W_{2}^{1,2}(Q)$.
\end{theorem}

If $T<\infty$ then the above theorem can be stated as follows.

\begin{theorem} \hspace{-5pt}.
\label{th7} Assume that  $T<\infty$, $a_i\in L_\infty(G)$ $(i=0,1,\ldots,n)$, and
the conditions  {\rm \eqref{eq601}-\eqref{eq603}}, where $\lambda=0$, hold.
 Then there exists  a unique solution to the problem
 \begin{equation}\label{eq602}
w_{0t}+Lw_0 =f_{0}(t,x), \ \   Bw_0|_{S}=g,\  \ w_0|_{t=0}=u_{0}(x),\
 \end{equation}
such that  $ w_{0}\in W_{2}^{1,2}(Q)$.
\end{theorem}

These theorems are consequences  of the results presented in  \cite[Theorems 5.7,7.11,8.2]{denk}.

Consider the problem  \eqref{eq21}-\eqref{eq23}. After the change of variables  $w=u-w_{0}$, we arrive at the simpler problem
 \begin{equation}\label{em3}
w_{t}+Lw = \sum_{i=1}^r q_i(t)\delta(x-x_i),
 \end{equation}
\begin{equation}\label{em4}
    Bw|_{S}=0,\  \ w|_{t=0}=0,\
    \end{equation}
\begin{equation}\label{em5}
    \ w(b_{j},t)=\psi_{j}(t)-w_{0}(t,b_{j})=\tilde{\psi}_{j}(t),\  j=1,2,\ldots,s.
    \end{equation}
Assume that the functions $\tilde{\psi}_{j}(t)$ admits the representation
  \begin{equation}\label{em6}
 \tilde{\psi}_{j}(t)=
   \int_{0}^{t}V_{\delta_{j}}(t-\tau)\psi_{0j}(\tau)d\tau,\     \psi_{0j}e^{-\lambda t}\in L_{2}(0,T),
    \end{equation}
    where $V_{\gamma}(t)=\frac{e^{-\gamma^{2}/4t}}{4\pi t}$ for $n=2$ and
 $V_{\gamma}=\frac{\gamma
e^{-\gamma^{2}/4t}}{2\sqrt{\pi}t^{3/2}}$ for $n=3$.

Assume that   $K=\{y\in G:\ \rho(y,\cup_{i=1}^{m}x_{i})\leq
\rho(y,\Gamma)\}$ for $G\neq {\mathbb R}^{n}$ and $K$ is an arbitrary compact set otherwise.
Denote  $$
\varphi_{j}(x)=\frac{-1}{2}\int_{0}^{1}(\vec{a}(b_{j}+\tau
(x-b_{j})),(x-b_{j}))\,d\tau.
$$
Let   $\delta_{j}=min_{i}r_{ij}, j=1,2,\ldots,s$, where
$r_{ij}=|x_{i}-b_{j}|$. Introduce the matrix   $A_{0}$ with the entries
$a_{ji}=e^{\varphi_{j}(x_{i})}$ if $|x_{i}-b_{j}|=\delta_{j}$
and   $a_{ji}=0$ otherwise.  We need the conditions
 \begin{equation}\label{ew16}
det{A}_{0}\neq 0.
 \end{equation}
\begin{equation}\label{ew17}
 \nabla \varphi_j, \Delta \varphi_j\in L_\infty(G), j=1,\ldots, r
\end{equation}
We also assume that all coefficients in \eqref{eq21} admits extensions
to the whole
 ${\mathbb R}^{n}$ such that the conditions \eqref{q13}, \eqref{ew17} hold for the case of  $G={\mathbb R}^{n}$.
If $G$ is a domain with compact boundary such an extension always exists.

Fix $p\in (1,n/(n-1))$. Let the space $W_{q,B}^{1}(G)$  agrees with $W_{q}^{1}(G)$ in the case of the
Robin boundary conditions and with the subspace of $W_{q}^{1}(G)$ comprising the functions vanishing on $\Gamma$ otherwise.
$W_{p,B}^{-1}(G)$ is the dual space to $W_{q,B}^{1}(G)$, $1/p+1/q=1$.

\begin{theorem} \hspace{-5pt}.
\label{th8} Assume that  $T=\infty$, $r=s$,
the conditions
 {\rm \eqref{eq601},  \eqref{q13}, \eqref{ew16}, \eqref{ew17}}  hold and  $b_{i}\in K$ for  $i=1,2,\ldots,s$.
 Then there exists  $ \lambda_{0}\geq 0$ such that if
 $\lambda\geq \lambda_{0}$ and the conditions
{\eqref{eq603}, \eqref{eq604}, \rm \eqref{em6}} are fulfilled then there exists a unique solution
to the problem
 {\rm \eqref{eq21}-\eqref{eq23}} such that  $u=w_{0}+w$, $w_{0}$ is a solution to
 the problem  {\rm \eqref{eq602}}, $e^{-\lambda t}w_{0}\in W_{2}^{1,2}(Q)$,
$e^{-\lambda t}\vec{q}\in L_{2}(0,\infty)$,  $e^{-\lambda t} w\in L_{2}(0,\infty;W_{p,B}^{1}(G))$,
$e^{-\lambda t} w_{t}\in L_{2} (0,\infty;W_{p,B}^{-1}(G))$,  $e^{-\lambda t} w\in W_{2}^{1,2}(Q_{\varepsilon})$  for all   $\varepsilon>0$.
\end{theorem}

If $T<\infty$ then Theorem \ref{th8} can be stated as follows.

\begin{theorem} \hspace{-5pt}.
\label{th9} Assume that  $T<\infty$, $r=s$,
the conditions
 {\rm \eqref{eq604}-\eqref{eq601}} with $\lambda=0$,  {\rm \eqref{q13}, \eqref{ew16}, \eqref{ew17}} hold and  $b_{i}\in K$ for  $i=1,2,\ldots,s$.
 Then there exists  a unique solution
to the problem
 {\rm \eqref{eq21}-\eqref{eq23}} such that  $u=w_{0}+w$, $w_{0}$ is a solution to
 the problem  {\rm \eqref{eq602}}, $w_{0}\in W_{2}^{1,2}(Q)$,
$\vec{q}\in L_{2}(0,\infty)$,  $ w\in L_{2}(0,\infty;W_{p,B}^{1}(G))$,
$w_{t}\in L_{2} (0,\infty;W_{p,B}^{-1}(G))$,  $ w\in W_{2}^{1,2}(Q_{\varepsilon})$  for all   $\varepsilon>0$.
\end{theorem}

The proof can be found in \cite[Theorem 2]{neus2}, \cite[Theorem 3]{neus4}.
It is based on the Laplace transform and asymptotic representations of a solution to the problem (see \cite[Theorem 7]{neus2}, \cite{pya112}).
\begin{equation}\label{eq24}
\lambda w+Lw = \delta(x-x_1),
 \end{equation}
\begin{equation}\label{eq25}
    Bw|_{\Gamma}=0
    \end{equation}
The scheme of the proof is presented in the Theorem \ref{th4} above.

We note that the following condition is actually a necessary
condition for the uniqueness of solutions to the problem
\eqref{eq1}-\eqref{eq3}. If it fails then any number of the
points $\{b_{i}\}$ does not ensure uniqueness of solutions
 (see examples below).

{Condition (D).} For $n=2$, any three points $\{b_{i}\}$ do not
lie on the same straight line and, for $n=3$, any four points
$\{b_{i}\}$ do not lie on the same plane.

{\bf An algorithm of determination of a solution   $u,q_{1}(t),x_{1}$ to the problem
\eqref{eq21}-\eqref{eq23} in the case of $r=1$.}
Let  $G={\mathbb  R}^{n}$ or  $G={\mathbb  R}^{n}_{+}$, or $G$ is a domain with compact boundary.
 We consider compact $K\subset G$ with the following
properties: $K$ contains the points $\{b_{i}\}_{i=1}^s$,  if $Bu=u$ and $G$ is a domain with compact boundary
then the convex hull of $K$ is  contained in $G$; if $G$ is a
domain with compact boundary and $Bu \neq u $ then $K\subset
B_{\rho(x_{1},\Gamma)}(x_{1})$; if $G={\mathbb R}^{n}$  then $K\subset G$ is an arbitrary compact set;
if $G={\mathbb R}^{n}_{+}$, $Bu=u$ or $Bu=-\frac{\partial u}{\partial n}$  and $a_i=0$ for $i=1,2,\ldots,n$
then $K\subset G$ is an arbitrary compact set.
If $G\neq {\mathbb R}^{n}$ then we can also take $K=\{y\in G:\ \rho(y,x_{1})\leq
\rho(y,\Gamma)\}$ as in Theorem \ref{th8}.

These conditions on the compact set $K$ imply that on the set  $K_\varepsilon=\{x\in K:\ |x-x_1|\geq \varepsilon\}$
 a solution to the problem \eqref{eq24}, \eqref{eq25} for a sufficiently large $\lambda\in {\mathbb R}$ admits
 the asymptotic representation (see Theorems 3.9, 3.11-3.13 in \cite{pya112})
  \begin{equation}\label{po13}
w(x)=\frac{1}{2\sqrt{2\pi|x-x_{1}|}\lambda^{1/4}} e^{-\sqrt{\lambda}|x-x_{1}|}\Bigl(1+ O\Bigl(\frac{1}{\sqrt{|\lambda|}}\Bigr)\Bigr),\ n=2;
    \end{equation}
 \begin{equation}\label{po14}
w(x)=\frac{1}{4\pi|x-x_{1}|} e^{-\sqrt{\lambda}|x-x_{1}|}\Bigl(1+ O\Bigl(\frac{1}{\sqrt{|\lambda|}}\Bigr)\Bigr),\ \ n=3;
    \end{equation}

Assume that
the functions  $u$, $q_{1}$ belongs to the class described in Theorem \ref{th8},  the data  $g,u_{0},f_{0}$ satisfy the conditions of this theorem for some sufficiently large parameter $\lambda\in {\mathbb R}$, $\lambda\geq \lambda_0$ with the parameter  $\lambda_0 $  that of Theorem \ref{th6}, and
 the condition  (D) is fulfilled.

 First we consider the case of  $n=3$. As before, we reduce our problem to the problem \eqref{em3}- \eqref{em5} of the form
    \begin{equation}\label{p36}
   w_{t}+Lw=q_{1}(t)\delta(x-x_{1}),\ \  w(x,0)=0,\ \  Bw|_{S}=0,
    \end{equation}
        \begin{equation}\label{p37}
    w(b_{i},t)=\psi_{i}-e^{\lambda t}v_{0}(b_{i},t)=\tilde{\psi}_{i}(t),\ i=1,2,\ldots,s.
    \end{equation}
 Using the  Laplace transform, we obtain that the function
 $\hat{w}={\cal L}(w)=\int_{0}^{\infty}e^{-\lambda
t}w(x,t)\,dt$ is a solution to the problem
    \begin{equation}\label{q24}
    \lambda \hat{w}+L\hat{w}=\hat{q}_{1}(\lambda)\delta(x-x_{1}),\ \   B\hat{w}|_{\Gamma}=0.
    \end{equation}
        \begin{equation}\label{q25}
    \hat{w}(b_{j})={\cal L}(\tilde{\psi}_{j})(\lambda),\ j=1,\ldots,s.
    \end{equation}
In view of  \eqref{p37}, we infer
      \begin{equation}\label{q26}
L(\tilde{\psi}_{j})(\lambda)=\hat{q}_{1}(\lambda) v_{0}(b_{j},\lambda), \ \ j=1,\ldots,s,
    \end{equation}
where $v_{0}$ is a solution to the problem  \eqref{eq24}, \eqref{eq25}.
The equalities \eqref{po13},  \eqref{po14} yield
      \begin{equation}\label{q28}
\frac{L(\tilde{\psi}_{j})(\lambda)}{L(\tilde{\psi}_{i})(\lambda)}=\frac{e^{-\sqrt{\lambda}|x_{1}-b_{j}|}{|x_{1}-b_{i}|}}
{e^{-\sqrt{\lambda}|x_{1}-b_{i}|}{|x_{1}-b_{j}|}(1+O(\frac{1}{\sqrt{\lambda}}))},\ \ i,j=1,2,3,4.
    \end{equation}
Denote
$G_{ij}={\cal L}(\tilde{\psi}_{j})(\lambda)/{\cal L}(\tilde{\psi}_{i})(\lambda)$.
Without loss of generality, we can assume that
$|O(\frac{1}{\sqrt{\lambda}})|\leq 1/2$ for all $i,j$ and
$\lambda\geq \lambda_{1}$ $(\lambda_{1}\geq \lambda_{0})$. Denote
$\alpha=\sqrt{\lambda}$. The equalities  \eqref{q28} imply that
      \begin{equation}\label{q29}
\frac{G_{ij}(\alpha^{2})}{G_{ij}((\alpha+1)^{2})}=e^{|x_{1}-b_{j}|-|x_{1}-b_{i}|}(1+O(\frac{1}{\alpha})).
    \end{equation}
Thus, we have
 \begin{equation}\label{q30}
\ln\frac{G_{ij}(\alpha^{2})}{G_{ij}((\alpha+1)^{2})}=\alpha_{j}-\alpha_{i}+ O(\frac{1}{\alpha}), \ \alpha_{j}=|x_1-b_j|.
    \end{equation}
It easy to see that there exists the limit
 \begin{equation}\label{q31}
\lim_{\alpha\to +\infty}\ln\frac{G_{ij}(\alpha^{2})}{G_{ij}((\alpha+1)^{2})}=\alpha_{j}-\alpha_{i}.
    \end{equation}
It is possible that  $\alpha_{1}=\alpha_{2}=\alpha_{3}=\alpha_{4}$,
in this case all limits in  \eqref{q31} are equal to zero. In this case, the point  $x_{1}$ is the center of the sphere
containing the points
$b_{1},b_{2},b_{3},b_{4}$, which is uniquely determined. Next, we have that
$\hat{q}_{1}=\frac{{\cal L}(\tilde{\psi}_{j})(\lambda)}{v_{0}(b_{j},\lambda)}$,
$q_{1}={\cal L}^{-1}\Bigl(\frac{{\cal L}(\tilde{\psi}_{j})(\lambda)}{v_{0}(b_{j},\lambda)}
\Bigr)$.
Assume that there is a nonzero limit $\alpha_{j_{0}}- \alpha_{i_{0}}$  among those in  \eqref{q31}. We can write out the approximate equality
\begin{equation}\label{q33}
\ln\frac{G_{ij}(\alpha^{2})}{G_{ij}((\alpha+1)^{2})}\approx \alpha_{j_{0}}-\alpha_{i_{0}},
    \end{equation}
 The equality  \eqref{q28} can be rewritten in the form
  \begin{equation}\label{q32}
\frac{\alpha_{i_{0}}}{\alpha_{j_{0}}}\approx
 \frac{L(\tilde{\psi}_{j_{0}})(\lambda)}{L(\tilde{\psi}_{i_{0}})(\lambda)}e^{\sqrt{\lambda}(\alpha_{j_{0}}-\alpha_{i_{0}})}
 =G(\lambda).
    \end{equation}
Sine the left-hand side is not identically  1, there exists
 $\lambda_{3}>0$ such that  $G(\lambda)\neq 1$ for  $\lambda\geq \lambda_{3}$. We obtain approximate  system   \eqref{q33},
   \eqref{q32} for determination of
the quantities  $\alpha_{j_{0}},\alpha_{i_{0}}$, and  the determinant of this system does not vanish uniformly on $\lambda$ such that
 $\lambda\geq \lambda_{4}$ for some sufficiently large  $\lambda_{4}$.
 A solution to  this system
 is some approximations of the numbers   $\tilde{\alpha}_{j_{0}},\tilde{\alpha}_{i_{0}}$. Hence, we have the representations
 $\alpha_{i_{0}}=\tilde{\alpha}_{j_{0}}+
O(\frac{1}{\alpha})$. Choosing
the nonzero limits  $\alpha_{i}-\alpha_{j}$ or using the equalities  $\alpha_{i}=\alpha_{j}$ in the case of this limit is zero, we can find
the quantities  $\alpha_{i}$, $i=1,2,3,4$.
In this case the point $x_{1}$ is the point of intersection of the spheres centered
at $b_{i}$ of radius  $\alpha_{i}$. The function $q_{1}$ can be determined with the use of the inverse  Laplace transfrom
 (or with the use of Duhamel formula), and the function  $w(x,t)$ is a solution to the problem \eqref{p36}, \eqref{p37}.

Consider the case of  $n=2$. Instead of  \eqref{q28}, we have
\begin{equation}\label{q34}
\frac{L(\tilde{\psi}_{j})(\lambda)}{L(\tilde{\psi}_{i})(\lambda)}=\frac{e^{-\sqrt{\lambda}|x_{1}-b_{j}|}{\sqrt{|x_{1}-b_{i}|}}}
{e^{-\sqrt{\lambda}|x_{1}-b_{i}|}{\sqrt{|x_{1}-b_{j}|}}(1+O(\frac{1}{\sqrt{\lambda}}))},\ \ i,j=1,2,3.
    \end{equation}
The equality  \eqref{q30} is of the form
 \begin{equation}\label{q35}
\ln\frac{G_{ij}(\alpha^{2})}{G_{ij}((\alpha+1)^{2})}=\alpha_{j}-\alpha_{i}+ O(\frac{1}{\alpha}).
    \end{equation}
The equality  \eqref{q32} can be written as
  \begin{equation}\label{q36}
\sqrt{\frac{\alpha_{i_{0}}}{\alpha_{j_{0}}}}=
 \frac{L(\tilde{\psi}_{j_{0}})(\lambda)}{L(\tilde{\psi}_{i_{0}})(\lambda)}e^{\sqrt{\lambda}(\alpha_{j_{0}}-\alpha_{i_{0}})}
 =G(\lambda).
    \end{equation}
The remaining arguments are similar.

\begin{remark} The numerical experiments show that these algorithms can be used for determination of the point $x_1$.
However, the procedure  requires an appropriate choice of the parameter $\lambda$ \cite{pya25}.
\end{remark}

Next, we describe some model situation in which $Lu=-\Delta u
+\lambda_{0}u$, $\lambda_{0}\geq 0$, $G={\mathbb R}^{n}$ and
functions $q_{i}$ on the right-hand side of \eqref{eq21} are
real constants.

\begin{theorem} \label{th10} Let $u_{1}, u_{2}$ be two solutions to the problem {\rm \eqref{eq21}-\eqref{eq23}} from
class described in the theorem {\rm \ref{th8}} with the
right-hand sides in \eqref{eq21} of the form
$\sum_{i=1}^{r_{j}}q_{i}^{j}\delta(x-x_{i})$ $(q_i^{j}=const,
j=1,2)$, the condition {\rm (D)} holds, and $s\geq 2r+1$ in the
case $n=2$ and $s\geq 3r+1$ in the case $n=3$, where $r\geq
\max(r_{1},r_{2})$ (i.e., there is the upper bound for   the
number $\max(r_{1},r_{2})$). Then $u_{1}=u_{2}$, $r_{1}=r_{2}$,
and $q_{i}^{1}=q_{i}^{2}$ for all $i$, i.~e.,   a solution to
the problem of recovering the number $m$, points $x_{i}$, and
constants $q_{i}$ is unique.
 \end{theorem}

We now  display  the corresponding examples showing the accuracy the
results obtained. The following example shows that if  the
condition (D) fails then the problem of recovering
 the intensities  of sources (sinks) located at $x_{1},x_{2}$
 has a nonunique solution. At the same time,
it is an example of the nonuniqueness in the problem of
recovering the intensity of one source and its location.
 Note that the problem of determining the location of one
source $x_{0}$ and its intensity $q(t)$ is simple enough  and to uniquely recover these parameters  we need two measurements in the case of $n=1$ \cite{pya1}, three measurements  in the case of $n=2$
\cite{pya25} and four measurements (that is $s=4$ in \eqref{eq3}) in the case of $n=3$ \cite{pya33}. The smaller number of points does not allow to define the parameters $q(t),x_{0}$ uniquely. We
should also require that  the point $x_{0}$ lie   between two measurement points in the case of $n=1$ and the condition (D) holds in the case of $n=2,3$. The numerical solution of the problem of
recovering one source is treated in the articles \cite{pya25},  \cite{18}-\cite{liu}.

\begin{example} First we take  $n=3$, $G={\mathbb R}^{n}$,
$Lu=-\Delta u$. Let $u$ be a  solution to the equation
\eqref{eq21} satisfying the homogeneous initial conditions  with
the right-hand side in $\eqref{eq21}$ of the form
$$
q(t)(\delta(x-x_{1})-\delta(x-x_{2})).
$$
The Laplace transform of this solution to the problem
\eqref{eq21}-\eqref{eq23} is written  as
$$
\hat{u}=\hat{q}(\lambda)(\frac{1}{4\pi|x-x_{1}|} e^{-\sqrt{\lambda}|x-x_{1}|}- \frac{1}{4\pi|x-x_{2}|} e^{-\sqrt{\lambda}|x-x_{2}|}).
$$
Let $P$ be the plane  perpendicular to the segment
$[x_{1},x_{2}]$ and passing through its center. We have
$$
\hat{u}(y,\lambda)\equiv 0 \ \ \forall y\in P.
$$
So,  $u(y,t)=0$ for all $y\in P$. Precisely the same example
 can be constructed  in the
case $n=2$. We take the  perpendicular to the segment
$[x_{1},x_{2}]$ passing through its center rather than  the
plane  $P$.  Thus, if condition (D) fails then any number of
measurement points does not allow to determine the intensity
and the location of the sources.
\end{example}

\begin{example} Consider the case of $G={\mathbb R}^{n}$,
$Lu=-\Delta u$. In this case the conditions \eqref{eq23}
with $s=4$ in the case of  $n=2$ and  $s=6$ in the case of
$n=3$ does not allow to determine location of two sources and
their intensities even if the condition (D) holds. We describe the example in the case of $n=3$. The case of $n=2$ is quite similar.  Let
$u_{1}$, $u_{2}$ be  solutions to the equation \eqref{eq21}
satisfying the homogeneous initial conditions in which the
right-hand sides are of the form
$$
q(t)\delta(x-x_{1})+q(t)\delta(x-x_{2}), \ q(t)\delta(x-x_{1}^{*})+q(t)\delta(x-x_{2}^{*}).
$$

The Laplace transforms of
$\hat{u}_{1}, \hat{u}_{2}$ are as follows:
  \begin{equation}\label{ew42}
\hat{u}_{1}(x,\lambda)=\sum_{i=1}^{2}\frac{\hat{q}}{4\pi|x-x_{i}|}
e^{-\sqrt{\lambda}|x-x_{i}|} , \ \ \hat{u}_{2}(x,\lambda)=\sum_{i=1}^{2}\frac{\hat{q}}{4\pi|x-x_{i}^{*}|}
e^{-\sqrt{\lambda}|x-x_{i}^{*}|}.
    \end{equation}
 Here we use explicit representations of the fundamental
solution for the Helmgoltz equation.
 We take   $x_{1}=(a,a,0), \  x_{1}^{*}=(a,-a,0), \
x_{2}=(-a,-a,0), \  x_{2}^{*}=(-a,a,0)$ ($a>0$). As is easily
seen,  the functions $\hat{u}_{1}, \hat{u}_{2}$ coincide at the
points
\begin{multline*}
b_{1}=(M,0,0),\  b_{2}=(-M,0,0),\  b_{3}=(0,M,0),\\ b_{4}=(0,-M,0), \ b_{5}=(0,0,M), \ b_{6}=(0,0,-M),
\end{multline*}
where $M>0$ and, thus, the problem of recovering the locations
of 2 sources and their intensities admits several solutions in
the case of $s=6$. It  follows from the theorem \ref{th10} that
 in the case of $s=7$ points $x_{1}, x_{2}$ and the intensities are
determined uniquely (if  the condition (D) holds  and  the
intensities are constants).
\end{example}

\vspace{5mm}
{ \bf Acknowledgement.}
{\it The research was carried out within the state assignment of Ministry of Science and Higher Education of the Russian Federation (theme No. FENG-2023-0004,
"Analytical and numerical study of inverse problems on recovering parameters of atmosphere or water pollution sources and (or) parameters of media")
 }

\vspace{5mm} \noindent {Sergey  Pyatkov, Doctor of Sci.,  Professor, Digital Technologies Department, Yugra State University,  Khanty-Mansiysk, Russia, E-mail: pyatkovsg@gmail.com \\ Lyubov Neustroeva, postgraduate, Digital Technologies Department, Yugra State University,  Khanty-Mansiysk, Russia, E-mail: starkovaLV@mail.ru}


\begin{thebibliography}{23}

\bibitem{mar}  G.I. Marchuk, {\it Mathematical Models in Environmental Problems}, %
    Studies in Mathematics and its Applications,
    Elsevier Science Publishers,  Amsterdam,{\bf 16}, 1986.

\bibitem{ozi}   M.N. Ozisik,   H.R.B. Orlande, {\it Inverse Heat
    Transfer,}  Taylor \& Francis, New York, 2000.

\bibitem{ali}     O.M. Alifanov, {\it  Inverse heat transfer  problems},
    Berlin, Springer-Verlag, 1994.

  \bibitem{pri}   A.I. Prilepko,  V.V. Solov’ev Solvability theorems and
{\it Rothe’s method for inverse Problems for a Parabolic Equation}, I,  Differ. Equations,  {\bf 23} (1987), 1230-1237.

    \bibitem{sol}  V.V. Solov’ev, {\it  Global Existence of a Solution to the Inverse
Problem of Determining the Source Term in a Quasilinear Equation of Parabolic Type},  Differ. Equations, {\bf  32}(1996),  538-547.

    \bibitem{pri2}    A.I. Prilepko,   D.G. Orlovsky,   I.A. Vasin {\it   Methods for solving inverse
        problems in     Mathematical Physics}.  Marcel Dekker, Inc.,  New-York,   1999.


\bibitem{bib3}   M. Ivanchov,  {\it Inverse Problems for Equation of
    Parabolic Type}, Math. Studies. Monograph Series, {\bf ~10},
  WNTL Publishers,     Lviv,  2003.

\bibitem{pya3}    S.G. Pyatkov,    V.V. Rotko, {\it   On
    some parabolic inverse problems with the pointwise
    overdetermination},   AIP Conference Proceedings,  {\bf 1907} (2017),  020008.

\bibitem{pya25}  Egor Safonov, Sergey Pyatkov,  {\it  On Some Classes of Inverse Problems on Determining the Source Function}. Proceedings of the 8th Scientific Conference on Information Technologies
    for     Intelligent Decision Making Support (ITIDS 2020), Series:  Advances in Intelligent Systems Research,
     Atlantis Press SARL,  Paris,  {\bf 483} (2020), 116-120.

\bibitem{pya6}  S.G. Pyatkov, V.A. Baranchuk,  {\it On some Inverse Parabolic Problems with Pointwise Overdetermination},  Journal of Siberian Federal University, Mathematics \& Physics,
    {\bf  14(4)} (2021),  463-474.
    % Журнал Сибирского федерального университета.Математика и физика.

\bibitem{pya7}   S.G. Pyatkov,  V.V. Rotko,  {\it Inverse problems with Pointwise Overdetermination for some quasilinear parabolic systems},  Siber. Adv. in Math., {\bf 30 (2)} (2020),
    124-142.

\bibitem{16}  V.V. Penenko, {\it Variational methods of data
    assimilation and inverse problems for
studying the atmosphere, ocean, and environment},  Numerical Analysis and Applications, {\bf 2(4)} (2009),  341-351.

\bibitem{18}  X. Deng,  Y. Zhao, J. Zou,  {\it On linear finite
    elements for simultaneously recovering source location and
    intensity}  Int. J. Numer. Anal. Model, {\bf 10} (2013), 588--602.

\bibitem{pen}  A.V. Penenko,  S. Rachmetullina, {\it Algorithms for
    atmospheric emission source localization based on the
    automated ecological monitoring system data},  Siberian
    Electronic Mathematical Reports,  10 (2013),  35-54.

\bibitem{4}  A.~El Badia,  T. Ha-Duong, ~A. Hamdi, {\it     Identification of a point source in a linear
    advection-dispersion-reaction
    equation: application to a pollution source problem},
Inverse Problems,  {\bf  21} (2005), 1121-1136.

 \bibitem{9}  A.~El Badia, ~A. Hamdi, {\it Inverse Source Problem in an Advection-Dispersion-Reaction System:
    Application to Water Pollution},  Inverse Problems, {\bf 23} (2007),  2103--2120.

 \bibitem{17} A.~El Badia, T. Ha-Duong,
{\it     Inverse source problem for the heat equation:
    application to a pollution detection problem}, J.~Inverse Ill-Posed Probl, {\bf 10} (2002),
 585--599.

\bibitem{bad}  A.~El Badia,  T. Ha-Duong, {\it  An inverse source
    problem in potential analysis}  Inverse Problems, {\bf 16} (2000), 651-663.

 \bibitem{11}  L. Ling,  T. Takeuchi,
{\it     Point sources identification problems for heat equations},  Commun. Comput. Phys., {\bf  5} (2009),  897--913.

\bibitem{pya1}  S.G. Pyatkov,  E.I. Safonov, {\it Point sources
    recovering problems for the one-dimensional heat equation},  Journal of Advanced Research in Dynamical and Control
    Systems,  {\bf 11(01)} (2019),  496-510.


\bibitem{neus}  L.V. Neustroeva,  {\it On uniqueness in the problems of determining point sources in mathematical models of heat and mass transfer},   Bulletin of the South Ural State University.
    Series Mathematics, Mechanics, Physics,  {\bf 14(2)} (2022),  31–43.

\bibitem{zhou}  L. Zhou,  P.K. Hopke,  W. Liu, {\it  Comparison of two
    trajectory based models for locating particle sources for
    two rural New York sites},  Atmospheric Environment,
{\bf 38} (2004),  1955-1963.


\bibitem{han}  Y. Han,  T.M. Holsen,  P.K. Hopke,  J. Cheong,
 H. S.Yi  Kim, {\it Identification of source location for
    atmospheric dry deposition of heavy metals during
    Yellow-Sand events in Seoul, Korea in 1998 using hybrid
    receptor models},  Atmospheric Environment,  {\bf 38} (2004),  5353-5361.


\bibitem{pek}  N.J. Pekney,  C.I. Davidson,  L. Zhow,  P.K. Hopke, {\it     Application of PSCF and CPF to PMF-modeled sources of PM2,5
    in Pittsburgh},  Aerosol Science and Technology,  {\bf  40} (2006),  952-961.

\bibitem{yan}  C.-Yu. Yang,  {\it Solving the two-dimensional heat
    source  problem through the linear  least square error
    method},  Int. J. Heat Mass Transfer, {\bf  41} (1998),  393-398.

\bibitem{mam}  A.V. Mamonov,   Y-H.R. Tsai,  {\it Point source
        identification in nonlinear advection-diffusion-reaction
        systems},  Inverse Problems,  {\bf 29(3))} (2013),  26.


\bibitem{ver}   N. Verdiere, G. Joly-Blanchard,   L. Denis-Vidal, {\it
    Identifiability and identification of a pollution source in
    a river by using a semi-discretized model},   Applied
    Mathematics and Computation, {\bf  221} (2013), 1-9.

\bibitem{maz}  M. Mazaheri,   J.M.V. Samani,  H.M.V. Samani, {\it
    Mathematical model for pollution source identification in
    rivers},  Environmental Forensics, {\bf  16(4)} (2015),     310-321.

\bibitem{su} J. Su , {\it Heat source estimation with the conjugate
    gradient method in inverse linear diffusive problems},
 J. Braz. Soc. Mech. Sci, {\bf 23(3))} (2001).

\bibitem{wang}  S. Wang,  L. Zhang,  X. Sun,  H.~Jial, {\it
 Solution to two-dimensional steady inverse heat transfer problems with interior heat source based on   the conjugate gradient method},
 Mathematical Problems     in Engineering, {\bf  2017} (2017),  2861342.

\bibitem{net}  A.J.S. Neto,  M.N. Ozisik,  {\it Two-dimensional
    inverse heat conduction problem of estimating the
    time varying strength of a line heat source},  Journal of
    Applied Physics, {\bf 71} (1992),  53-57.

\bibitem{mil}  E. Milnes,  {\it Simultaneous identification of a single
    pollution point-source location and contamination time
    under known flow field conditions},  Advances in Water
    Resources, {\bf  30} (2007), 2439-2446.

\bibitem{liu}  F.B. Liu,  {\it A modified genetic algorithm for solving
    the inverse heat transfer problem of estimating plan heat
    source},  International Journal of Heat and Mass Transfer,  {\bf  51(15-16)} (2008),  3745-3752.


\bibitem{hon} Y.C. Hon, M. Li, Y.A. Melnikov, {\it Inverse source identification by Green's function},
  Engineering Analysis with Boundary Elements,
{\bf  34(4)}  (2010),  352-358.

\bibitem{pya111}
 S.G. Pyatkov,  E.I. Safonov,  {\it On some classes of inverse problems of recovering a source function},  Siberian Advances in Mathematics,  {\bf 27(2)}  (2017), 119–132.


\bibitem{pya33}  L.V. Neustroeva, S.G. Pyatkov,   {\it On
    recovering a point source in some heat and
    mass transfer problems}, { AIP Conference Proceedings},
    {\bf 2328}  (2021), 020006.

\bibitem{neus1} L. Neustroeva, S. Pyatkov, {\it On solvability of some inverse problems of recovering point sources},  AIP Conference Proceedings, {\bf 2528}  (2022), 020005.

\bibitem{neus2}  S.G. Pyatkov,  L.V. Neustroeva,  {\it  On solvability of inverse problems of recovering pointwise sources}, Mathematical Notes of NEFU, {\bf  29(2)}  (2022),  43-58.

\bibitem{neus3} S.G. Pyatkov,  L.V. Neustroeva,  {\it  On some classes of inverse problems of recovering a source function}, Mathematical Notes of NEFU, {\bf  27(1)}  (2020),  21-40.

\bibitem{neus4}   L.V. Neustroeva, {\it  On uniqueness in the problems of determining point sources in mathematical models of heat and mass transfer},
 Bulletin of the South Ural State University. Series Mathematics, Physics, Mechanics, {\bf   14(2)}  (2022), 31–43.


\bibitem{pya112}  S.G. Pyatkov,  L.V. Neustroeva, {\it  On some asymptotic representations of solutions to elliptic equations and their applications}, Complex Variables and Elliptic Equations ,
  {\bf   66(6-7)}  ( 2021),     964-987.

\bibitem{ama}  H. Amann, {\it  Compact embeddings of vector-valued Sobolev and Besov spaces}, {\it Glasnik Mat.}, Ser.\,{\rm III}, {\bf  35(55)} (2000), 161-177.

  \bibitem{lad}    O.A. Ladyzhenskaya,  V.A. Solonnikov,  N.N. Ural’tseva,  {\it Linear and Quasilinear Equations of
Parabolic Type},   American Math. Society, Providence, R.I., 1968.

 \bibitem{agr}
 M.S. Agranovich, M.I. Vishik,
{\it     Elliptic problems with a parameter and parabolic problems of general type},
 Russ. Math. Surv. {\bf  19(3)} (1964), 53-157.

 \bibitem{denk}
 R. Denk,  M. Hieber,  J. Pr\"{u}ss,
{\it  R-boundedness, Fourier multipliers and problems  of elliptic and parabolic type},
   Mem. Am. Math. Soc.,   {\bf 166} (2003).


\end{thebibliography}
\end{document}